\documentclass[11pt]{article}
\usepackage[a4paper]{anysize}\marginsize{3.5cm}{3.5cm}{1.3cm}{2cm}
\pdfpagewidth=\paperwidth \pdfpageheight=\paperheight
\usepackage{amsfonts,amssymb,amsthm,amsmath,eucal}
\usepackage{pgf}
\usepackage{bbm}

\usepackage{tikz} 
\usepackage{float}
\usetikzlibrary{arrows}
\usepackage{subfigure}
\usepackage{caption}

\theoremstyle{plain}
\newtheorem{thm}{Theorem}[section]
\newtheorem{theorem}[thm]{Theorem}

\newtheorem{proposition}[thm]{Proposition}

\newtheorem{conjecture}[thm]{Conjecture}

\theoremstyle{definition}
\newtheorem{definition}[thm]{Definition}
\newtheorem{remark}[thm]{Remark}

\newtheorem{cexample}[thm]{Counterexample}

\newtheorem{question}[thm]{Question}
\newtheorem{problem}[thm]{Problem}

\newtheorem{thevarthm}[thm]{\varthmname}

\newenvironment{varthm*}[1]{\trivlist\item[]{\bf #1.}\it}{\endtrivlist}

\renewcommand\geq{\geqslant}

\renewcommand\leq{\leqslant}

\newcommand\be{\begin{eqnarray*}}
\newcommand\ee{\end{eqnarray*}}

\newcommand\R{\mathbb R}
\newcommand\C{\mathbb C}

\newcommand\K{\mathbb K}
\newcommand\F{\mathbb F}
\newcommand\LL{\mathbb L}
\renewcommand\P{\mathbb P}

\newcommand\cali{{\mathcal I}}

\newcommand\call{{\mathcal L}}
\newcommand\eps{\varepsilon}
\newcommand\newop[2]{\def#1{\mathop{\rm #2}\nolimits}}
\newop\log{log}
\newop\ord{ord}
\newop\Gal{Gal}
\newop\SL{SL}
\newop\Bl{Bl}
\newop\mult{mult}
\newop\mass{mass}
\newop\div{div}
\newop\codim{codim}
\newop\sing{sing}
\newop\Zeroes{Zeroes}
\newop\Ass{Ass}
\newop\reg{reg}
\newop\Ext{Ext}
\newop\PSL{PSL}
\newop{\symassreg}{areg}

\newcommand\lra\longrightarrow

\def\keywordname{{\bfseries Keywords}}%
\def\keywords#1{\par\addvspace\medskipamount{\rightskip=0pt plus1cm
\def\and{\ifhmode\unskip\nobreak\fi\ $\cdot$
}\noindent\keywordname\enspace\ignorespaces#1\par}}
\def\subclassname{{\bfseries Mathematics Subject Classification
(2000)}\enspace}
\def\subclass#1{\par\addvspace\medskipamount{\rightskip=0pt plus1cm
\def\and{\ifhmode\unskip\nobreak\fi\ $\cdot$
}\noindent\subclassname\ignorespaces#1\par}}

\definecolor{uuuuuu}{rgb}{0,0,1}
\definecolor{qqqqff}{rgb}{0,0,1}
\definecolor{xdxdff}{rgb}{0,0,1}

\def\endproof{\hspace*{\fill}\endproofsymbol\endtrivlist}

\def\endproofsymbol{\frame{\rule[0pt]{0pt}{6pt}\rule[0pt]{6pt}{0pt}}}

\begin{document}

\author{T.~Szemberg, J.~Szpond}
\title{On the containment problem}
\date{\today}
\maketitle
\thispagestyle{empty}

\begin{abstract}
   The purpose of this note is to provide an overview of the containment
   problem for symbolic and ordinary powers of homogeneous ideals,
   related conjectures and examples. We focus here on ideals with
   zero dimensional support. This is an area of ongoing active research.
   We conclude the note with a list of potential promising paths
   of further research.
\subclass{MSC 14C20 \and MSC 14J26 \and MSC 14N20 \and MSC 13A15 \and MSC 13F20}
\end{abstract}


\section{Introduction}
   Let $I\subset R=\K[x_0,\ldots,x_N]$ be a homogeneous ideal.
   Such an ideal determines various sequences of associated ideals.
   The most natural one is that of \emph{ordinary powers} of $I$:
   $$ I^0=R\supset I^1=I\supset I^2\supset I^3\supset I^4\supset \ldots.$$
   If $I=\langle f_1,\ldots,f_n\rangle$ is given explicitly in terms of generators,
   then the $r$-th ordinary power of $I$ is generated by products
   $$f_{i_1}\cdot\ldots\cdot f_{i_r}\;\mbox{ with }\; i_1,\ldots,i_r\in\left\{1,\ldots,n\right\}.$$
   Another important sequence of ideals associated to $I$ is the sequence of \emph{symbolic powers} of $I$:
   $$ I^{(0)}=R\supset I^{(1)}=I\supset I^{(2)}\supset I^{(3)}\supset I^{(4)}\supset \ldots.$$
   Recall that for $m\geq 0$ the $m$-th symbolic power of $I$ is defined as
   \begin{equation}\label{eq:symbolic power}
      I^{(m)}=\bigcap_{P\in\Ass(I)}\left(I^mR_P\cap R\right),
   \end{equation}
   where $\Ass(I)$ is the set of associated primes of $I$.

   Contrary to ordinary powers, even if $I$ is given explicitly by a set of generators
   $I=\langle f_1,\ldots,f_n\rangle$, it is very difficult in general to determine generators
   of $I^{(m)}$ for $m\geq 2$. On the other hand, if $I$ is a radical ideal
   and $V=V(I)$ is the reduced subscheme in $\P^N$ defined by $I$, then
   the symbolic power $I^{(m)}$ consists of all polynomials in $R$ vanishing
   to order at least $m$ along $V$. This is a consequence of the famous
   Nagata-Zariski Theorem, see \cite[Theorem 3.14]{Eisenbud} for prime ideals
   and \cite[Corollary 2.9]{SidSul09} for radical ideals. Thus, in this case,
   symbolic powers have clear geometric interpretation.

   Research presented in this note is motivated by the following main question.
\begin{problem}[The containment problem]\label{pro:containment}
   Decide for which $m$ and $r$ there is the containment
   \begin{equation}\label{eq:containment m in r}
      I^{(m)}\subset I^r.
   \end{equation}
\end{problem}
\begin{remark}
   Note that for the reverse containment i.e.
   $$I^r\subset I^{(m)}$$
   the analogous question has a very simple answer. The containment holds
   if and only if $m\leq r$ \cite[Lemma 8.4.1]{PSC}.
\end{remark}
   It came in a sense as a surprise that there does exist a uniform and simple
   stated answer to the Containment Problem. It was discovered by
   Ein, Lazarsfeld and Smith in characteristic zero, see \cite{ELS01}
   and by Hochster and Huneke in positive characteristic, see \cite{HoHu02}.
   Both papers builded upon ground breaking ideas of Swanson in \cite{Swa00}.
\begin{theorem}[Ein-Lazarsfeld-Smith, Hochster-Huneke]\label{thm:ELS}
   Let $I\subset\K[x_0,\ldots,x_N]$ be a homogeneous ideal. Then the containment
   $$I^{(m)}\subset I^r$$
   holds for all $m\geq Nr$.
\end{theorem}
   The above Theorem is a powerful and elegant result. It is natural to wonder
   to what extend the bound $m\geq Nr$ is sharp. A striking fact is that there is not
   a single example known (at least to the authors of the present note),
   where one really needs $m\geq Nr$ for the containment \eqref{eq:containment m in r}
   to hold, i.e. the containment holds for lower values of $m$.
   This has prompted Huneke to ask the following question addressing the first non-trivial
   case of Theorem \ref{thm:ELS}.
\begin{problem}[Huneke]\label{pro:Huneke}
   Let $I$ be a saturated ideal of a reduced finite set of points in $\P^2$.
   Does then the containment
   \begin{equation}\label{eq:containment 3 in 2}
      I^{(3)}\subset I^2
   \end{equation}
   hold?
\end{problem}
   Positive answer to Problem \ref{pro:Huneke} obtained in an ample family of special cases
   and additional experimental data led to the following more general question,
   see \cite[Conjecture 8.4.2]{PSC}, \cite[Conjecture 1.1]{BocHarPAMS10},
   \cite[Conjecture 4.1.1]{HaHu13}, \cite[Conjecture 3.2]{BCH14}.
\begin{problem}[Bocci, Harbourne, Huneke]\label{pro:BHH}
   Let $I$ be a saturated ideal of a finite set of reduced points in $\P^N$.
   Does the containment
   $$I^{(m)}\subset I^r$$
   hold for $m\geq Nr-(N-1)$?
\end{problem}
   This note gives an overview of verified cases and constructed
   counterexamples which have been discovered recently. There is also a number
   of open questions in this area of active current interest and ongoing investigations,
   see e.g. \cite{LBM15}, \cite{KKM15}, \cite{NagSec15} for recent contributions.
\section{Symbolic powers for ideals of points in projective spaces}
   Since in this note we focus on ideals supported on reduced or, sometimes more generally,
   on fat points, we recall here briefly the main notions. In the particular
   situation of points, symbolic powers can be defined in an easier and more
   accessible way than provided in \eqref{eq:symbolic power}.
   Let $Z=\left\{P_1,\ldots,P_s\right\}$ be a union of reduced points in $\P^N$. 
   For a fixed $i\in\left\{1,\ldots,s\right\}$
   let $I(P_i)$ denote the ideal of polynomials vanishing at $P_i$.
   Then obviously the ideal of the set $Z$ satisfies:
   $$I(Z)=I(P_1)\cap\ldots\cap I(P_s).$$
   The symbolic powers of $I(Z)$ are given in the same way:
   \begin{equation}\label{eq:symbolic Z}
      I(Z)^{(m)}=I(P_1)^{m}\cap\ldots\cap I(P_s)^{m}
   \end{equation}
   for any $m\geq 1$.
   Note that any point $P\in\P^N$ is a complete intersection,
   hence $I(P)^{(m)}=I(P)^m$ for all $m\geq 1$.
   Indeed since the powers of a complete intersection ideal are arithmetically
   Cohen-Macaulay the claim follows from the 
   Unmixedness Theorem \cite[Corollary 18.14]{Eisenbud}.
   So that the right hand side of \eqref{eq:symbolic Z}
   is the intersection of symbolic powers of ideals of each
   component of $Z$.
   
   The same holds in a slightly more general setting.
\begin{definition}[Fat points scheme]
   A subscheme $Z\subset \P^N$ is a \emph{fat points scheme}
   if its ideal $I(Z)$ has the form
   $$I(Z)=I(P_1)^{m_1}\cap\ldots\cap I(P_s)^{m_s}$$
   for points $P_1,\ldots,P_s$ in the projective space $\P^N$
   and positive integers $m_1,\ldots,m_s$.
   
   In this situation $I(Z)$ is a called a \emph{fat points ideal}.
\end{definition}   
   For symbolic powers of a fat points ideal $I(Z)=I(P_1)^{m_1}\cap\ldots\cap I(P_s)^{m_s}$ we have
   $$I(Z)^{(m)}=I(P_1)^{mm_1}\cap\ldots\cap I(P_s)^{mm_s}.$$
\section{General statements motivated by Problem \ref{pro:BHH}}
   To begin with note that the constant $N$ appearing in the statement of
   Theorem \ref{thm:ELS} can not be lowered. A series of examples was constructed
   by Bocci and Harbourne in \cite{BocHarJAG10}. The main idea to study
   the so called star configurations of points in $\P^N$ goes back to Ein.
\begin{definition}[Star configuration of points]\label{def:star config}
   We say that $Z\subset\P^N$ is a \emph{star configuration} of degree $d$
   (or a $d$-star for short)
   if $Z$ consists of \textbf{all} intersection points
   of $d\geq N$ \textbf{general} hyperplanes in $\P^N$. By intersection points
   we mean the points which belong to exactly $N$ of given $d$ hyperplanes.
\end{definition}
   The assumption \emph{general} in the Definition means that
   any $N$ of $d$ given hyperplanes meet in a single point
   and there is no point belonging to $N+1$ or more hyperplanes.
   Star configurations can be defined much more generally.
   They form an interesting and combinatorially easy to describe
   class of examples occurring in various situations in algebra
   and geometry. We refer
   to \cite{GHM13} for a very nice introduction to this circle of ideas.

   Now the result of Bocci and Harbourne is as follows, see \cite[Theorem 2.4.3]{BocHarJAG10}.
\begin{proposition}\label{prop:rho for d-star in PN}
   Let $Z\subset \P^N$ be a $d$-star.
   Let $I$ be the ideal of $Z$. Then for any $c<\frac{d-N+1}{d}N$
   there exists $m$ and $r$ such that
   $$\frac{m}{r}\leq c\;\mbox{ and }\; I^{(m)}\nsubseteq I^r.$$
\end{proposition}
   Since $d$ can be taken arbitrarily large, we see that there is
   no constant lower than $N$ which would satisfy the statement
   of Theorem \ref{thm:ELS}. Investigating into this problem has prompted
   Harbourne to introduce a new interesting invariant, measuring in
   a sense a discrepancy between symbolic and ordinary powers.
\begin{definition}[Resurgence]\label{def:resurgence}
   Let $I$ be a homogeneous ideal. The \emph{resurgence} of $I$ is the real number
   $$\rho(I)=\sup\left\{\frac{m}{r}:\; I^{(m)}\nsubseteq I^r\right\}.$$
\end{definition}
   This is a delicate invariant and has been computed only in few special cases.
   One of them is
   $$\rho(I)=\frac{d-N+1}{d}$$
   for a $d$-star in $\P^N$, see \cite[Theorem 2.4.3]{BocHarJAG10}.
   The importance of $\rho(I)$ follows from the fact that there is the containment
   \begin{equation}\label{eq:containment rho}
      I^{(m)}\subset I^r\;\mbox{ provided}\; \frac{m}{r}>\rho(I).
   \end{equation}
   Of course the applicability of \eqref{eq:containment rho} is limited
   by the computability of $\rho(I)$. There are however results along these lines
   which in certain situations lead to interesting consequences.
   Before stating sample results of this kind we need to introduce further
   invariants.
\begin{definition}[The initial degree]
   Let $I=\oplus_{d\geq 0}I_d$ be a homogeneous ideal in $R$. Then the
   \emph{initial degree} $\alpha(I)$ is defined as the least number $d$
   such that $I_d\neq 0$. In other words, $\alpha(I)$ is the least degree
   of a non-zero polynomial in $I$.
\end{definition}
   Alternatively, denoting by $M=\langle x_0,\ldots,x_N\rangle$ the irrelevant ideal of $R$,
   the initial degree $\alpha(I)$ may thought of as the $M$-adic order of $I$, i.e.
   the largest $t$ such that $I\subset M^t$.

   The next invariant, the Castelnuovo-Mumford regularity $\reg(I)$ of $I$
   provides a measure, in terms of the minimal free resolution, of how complicated the ideal $I$ is.
\begin{definition}[Castelnuovo-Mumford regularity]
   Let $I\subset R$ be a homogeneous ideal and let
   $$0\lra\ldots\lra F_j\lra\ldots\lra F_0\lra I\lra 0$$
   be the minimal free resolution of $I$ over $R$.
   Let $f_j$ be the maximal degree of a generator in 
   a minimal set of generators of $F_j$. Then
   $$\reg(I)=\max\left\{f_j-j,\;j\geq 0\right\}.$$
\end{definition}
   The following useful result is Lemma 2.3.4 in \cite{BocHarJAG10} and Lemma 2.1 in \cite{BCH14}.
\begin{proposition}[Postulation Containment Criterion]\label{prop:postulation containment criterion}
   Let $I$ be a homogeneous (not necessarily saturated) ideal defining
   a $0$-dimensional subscheme in $\P^N$. If the inequality
   $$r\cdot \reg(I)\leq \alpha(I^{(m)})$$
   holds, then there is the containment
   $$I^{(m)}\subset I^r.$$
\end{proposition}
\begin{remark}[Computability of $\alpha(I)$]
   Note that the initial degree of a symbolic power of an ideal
   is not easy to compute in general. Suffices it to say that
   the famous Nagata Conjecture predicts
   $$\alpha(I^{(m)})>{m}\cdot {\sqrt{s}}$$
   for $s\geq 10$ general points in $\P^2$.
   Despite tremendous efforts this conjecture remains open for more than
   half a century, see \cite{CHMR13} for a nice overview and recent progress.
\end{remark}
   Turning back to Problem \ref{pro:BHH} we note that the constant $N-1$ appearing there
   is the most optimistic one, since again star configurations can be used in order to
   show that the statement would fail with $N$ in place of $N-1$. On the other hand
   there is number of cases where the question asked in Problem \ref{pro:BHH}
   has positive answer.
\begin{theorem}[Evidence for Problem \ref{pro:BHH}]\label{thm:evidence}
   The containment in Problem \ref{pro:BHH} holds
   \begin{itemize}
      \item[a)] for arbitrary ideals in characteristic $2$;
      \item[b)] for monomial ideals in arbitrary characteristic;
      \item[c)] for ideals of $d$-stars;
      \item[d)] for ideals of general points in $\P^2$ and $\P^3$.
   \end{itemize}
\end{theorem}
\proof
   Claim a) follows from a private communication from Huneke to Harbourne
   mentioned in \cite[Example 8.4.4]{PSC}. Claim b) is an application
   of the Postulation Containment Criterion as explained
   in \cite[Example 8.4.5]{PSC}. Part c) is again an application of
   Proposition \ref{prop:postulation containment criterion}
   based on numerical data provided in \cite[Lemma 8.4.7]{PSC}.
   Finally claim d) is proved for points in $\P^2$ in \cite[Proposition 6.10]{HaHu13}
   (it follows also from \cite[Remark 4.3]{BocHarJAG10})
   and for points in $\P^3$ in \cite[Theorem 3]{Dum15}.
\endproof
   The first counterexample to Problem \ref{pro:Huneke} and hence also to Problem
   \ref{pro:BHH} was announced in \cite{DST13b}. We report on it in details
   in the subsequent section. Now we mention counterexamples to
   Problem \ref{pro:BHH} announced recently by Harbourne and Seceleanu in \cite{HarSec15}.
   These are up to date the only counterexamples known in higher dimensional
   projective spaces. Such examples are available only in positive characteristic.

   From now on let $\K$ be a field of odd characteristic $p$ and let $\LL$ be its
   subfield of order $p$.
\begin{cexample}\label{cex:finite char in I2}
   Let $N=\frac{p+1}{2}$ and let $Z$ be the set of all but one $\LL$-points in $\P^N(\K)$.
   Then for the ideal $I=I(Z)$ there is
   $$I^{(\frac{p+3}{2})}\nsubseteq I^2.$$
\end{cexample}
\proof
   The key observation is that whereas $\alpha(I^2)=p^2+1$, there is
   in the ideal $I^{(\frac{p+3}{2})}$ a form of degree $p^2$, see \cite[Theorem 3.9]{HarSec15}
   for details.
\endproof
   The second counterexample allows more flexibility in the choice of the ordinary power.
\begin{cexample}
   We choose now the numbers $p$ and $N$ so that $p\equiv 1\;(\mod N)$ and $p>(N-1)^2$.
   Let again $Z$ be the set of all but one $\LL$-points in $\P^N(\K)$.
   Then for $r=\frac{p-1}{N}+1$ there is
   $$I^{(p)}\nsubseteq I^r.$$
\end{cexample}
\proof
   See \cite[Theorem 3.10]{HarSec15}.
\endproof
   We conclude this section recalling a homological containment criterion
   invented recently by Seceleanu, see \cite[Theorem 3.1]{Sec15}.
   The criterion involves local homology.
   The underlying idea is quite simple,
   the containment $I^{(m)}\subset I^r$ induces a map between
   local cohomology modules $H^0_M(R/I^m)\to H^0_M(R/\cali^r)$.
\begin{theorem}[Seceleanu]\label{thm:Seceleanu}
   Let $I\subset R$ be a homogeneous ideal. Consider the associated exact sequence
   $$0\to I^r/I^m\to R/I^m\stackrel{\pi}{\to} R/I^r \to 0.$$
   Then the following conditions are equivalent:
   \begin{itemize}
   \item[i)] there is the containment $I^{(m)}\subset I^r$,
   \item[ii)] the induced map $H^0_M(\pi): H^0_M(R/I^m)\to H^0_M(R/I^r)$ is the zero map.
   \end{itemize}
\end{theorem}
   There is also the dual version involving the $\Ext$ functor. In fact this dual approach
   together with the minimal free resolutions of the involved ideals open door to an effective
   application of Theorem \ref{thm:Seceleanu}. Even in the case of points in a plane this
   is a non-trivial operation. Seceleanu developed additional arguments
   \cite[Theorem 3.3]{Sec15} in order to deal
   just with the case $\cali^{(3)}\subset \cali^2$.
   Her method was successfully applied in order
   to verify the non-containment for the Klein configuration of points and for the
   series of Fermat examples, see Section \ref{sec:cex to I3 in I2}.   
\section{Specific statements related to Problem \ref{pro:Huneke}}\label{sec:cex to I3 in I2}
   The first counterexample to Problem \ref{pro:Huneke}
   was announced in 2013 by Dumnicki, Tutaj-Gasi\'nska
   and the first author in \cite{DST13b}. Afterwards whole
   series of further counterexamples in all characteristics but $2$
   (see Theorem \ref{thm:evidence} a)) have been found.
   Strangely enough all these counterexamples arise taking as the support
   of $I$ all (or almost all) intersection
   points of configurations of lines either defined by reflection groups
   or distinguished in some other way.
   Such configurations are extremal also from the combinatorial point of view,
   e.g. they exhibit an unusually high number of intersection points of high
   multiplicity. This is a surprising phenomenon not understood yet.
\begin{cexample}[Dual Hesse configuration]\label{cex:dual hesse}
   Let $\eps$ be a primitive root of $1$ of order $3$.
   We consider the saturated and radical ideal $I$
   of the following set of $12$ points in the complex projective plane $\P^2$:
   \begin{equation*}
   \begin{array}{lll}
      P_1=(1:0:0),         & P_2=(0:1:0),         & P_3=(0:0:1),\\
      P_4=(1:1:1),         & P_5=(1:\eps:\eps^2), & P_6=(1:\eps^2:\eps),\\
      P_7=(\eps:1:1),      & P_8=(1:\eps:1),      & P_9=(1:1:\eps),\\
      P_{10}=(\eps^2:1:1), & P_{11}=(1:\eps^2:1), & P_{12}=(1:1:\eps^2).
   \end{array}
   \end{equation*}
   Then
   $$I^{(3)}\nsubseteq I^2.$$
\end{cexample}
\proof
   See \cite[Theorem 2.2]{DST13b}.
\endproof
   The points listed in Counterexample \ref{cex:dual hesse} are all
   intersection points of an arrangement of $9$ lines given explicitly
   by equations
   \begin{equation*}
   \begin{array}{lll}
   L_1:\; x-y=0,        & L_2:\; y-z=0,        & L_3:\; z-x=0,\\
   L_4:\; x-\eps y=0,   & L_5:\; y-\eps z=0,   & L_6:\; z-\eps x=0,   \\
   L_7:\; x-\eps^2 y=0, & L_8:\; y-\eps^2 z=0, & L_9:\; z-\eps^2 x=0.
   \end{array}
   \end{equation*}
   The product of these $9$ linear equations is an element in $I^{(3)}$ (since
   there are exactly $3$ lines passing through any of the points $P_1,\ldots,P_{12}$
   but this product is not contained in $I^2$. Up to projective change of coordinates
   there is just one configuration of this kind, see \cite{ArtDol09} for
   this and much more. It is also well known that such a configuration
   cannot exist in the real projective plane because it would violate the celebrated
   Sylvester-Gallai Theorem. Nevertheless a real counterexample has been discovered
   shortly after announcement of \cite{DST13b}. It is presented in \cite{Real}. 
   Quite expectedly this counterexample is related
   to constructions of real line configurations with the maximal number
   of triple points.
\begin{cexample}[B\"or\"oczky configurations]
   It is convenient to identify
   the real plane $\R^2$ with the set of complex numbers $\C$
   in the usual way. Let $n$ be an even and positive integer.
   Let $\xi=\exp(2\pi i/n)$ be a primitive $n$--th root of unity
   and let $P_i=\xi^{i}$ for $i=0,\ldots, n-1$ be vertices of a regular $n$--gon.
   For $m,k\in\left\{0,\ldots,n-1\right\}$ we denote by $L_{m,k}$ the (real) line
   determined by the points $P_{m}$ and $P_{k}$ if $m\neq k$
   and the tangent line to the unit circle at the point $P_{m}$ if $m=k$.

   The configuration of lines we are interested in is defined as the union of $n$ lines
   $$\call_n=\left\{L_{i,\frac{n}{2}-2i},\; i=0,\ldots n-1\right\},$$
   where the indices are understood $\mod n$.
   Let $Z_{n}$ be the set of all triple points in the configuration $\call_n$.
   Such a configuration, for $n=12$ is visualized in Figure \ref{fig:B12}.
   The $19$ points in the set $Z_{12}$ are marked by dots.
\begin{figure}[H]
\centering
\definecolor{uuuuuu}{rgb}{1,1,1}
\begin{tikzpicture}[line cap=round,line join=round,x=1.0cm,y=1.0cm]
\clip(-0.82,-3.96) rectangle (8.3,5.22);
\draw [domain=-4.36:18.32] plot(\x,{(-2.33-0.03*\x)/-2.89});
\draw [domain=-4.36:18.32] plot(\x,{(-5.21--0.51*\x)/-1.98});
\draw [domain=-4.36:18.32] plot(\x,{(--7.45-1.43*\x)/1.46});
\draw [domain=-4.36:18.32] plot(\x,{(--8.44-2.49*\x)/1.47});
\draw [domain=-4.36:18.32] plot(\x,{(--4.08-1.97*\x)/0.55});
\draw [domain=-4.36:18.32] plot(\x,{(-3.23--1.98*\x)/0.51});
\draw [domain=-4.36:18.32] plot(\x,{(-6.11--2.51*\x)/1.42});
\draw [domain=-4.36:18.32] plot(\x,{(-5.12--1.46*\x)/1.43});
\draw [domain=-4.36:18.32] plot(\x,{(--2.03-0.55*\x)/-1.97});
\draw [domain=-4.36:18.32] plot(\x,{(--5.19-0.71*\x)/1.26});
\draw [domain=-4.36:18.32] plot(\x,{(--3.18-0.73*\x)/-1.24});
\draw [domain=1:2] plot(\x,{(-2.12--1.44*\x)/-0.01});
\begin{scriptsize}
\fill [color=uuuuuu] (2.19,-0.42) circle (2.0pt);
\fill [color=uuuuuu] (4.35,0.85) circle (2.0pt);
\fill [color=uuuuuu] (2.17,2.08) circle (2.0pt);
\fill [color=uuuuuu] (0.9,4.23) circle (2.0pt);
\fill [color=uuuuuu] (1.44,3.32) circle (2.0pt);
\fill [color=uuuuuu] (1.45,2.26) circle (2.0pt);
\fill [color=uuuuuu] (2.36,2.8) circle (2.0pt);
\fill [color=uuuuuu] (1.85,0.82) circle (2.0pt);
\fill [color=uuuuuu] (2.9,0.83) circle (2.0pt);
\fill [color=uuuuuu] (1.48,-1.68) circle (2.0pt);
\fill [color=uuuuuu] (0.96,-2.6) circle (2.0pt);
\fill [color=uuuuuu] (2.39,-1.14) circle (2.0pt);
\fill [color=uuuuuu] (1.47,-0.62) circle (2.0pt);
\fill [color=uuuuuu] (3.44,-0.08) circle (2.0pt);
\fill [color=uuuuuu] (4.88,0.32) circle (2.0pt);
\fill [color=uuuuuu] (4.87,1.38) circle (2.0pt);
\fill [color=uuuuuu] (5.79,0.86) circle (2.0pt);
\fill [color=uuuuuu] (6.85,0.87) circle (2.0pt);
\fill [color=uuuuuu] (3.42,1.75) circle (2.0pt);
\fill [color=uuuuuu] (3.61,2.09) circle (2.0pt);
\fill [color=uuuuuu] (3.64,-0.41) circle (2.0pt);
\fill [color=uuuuuu] (1.46,0.82) circle (2.0pt);
\end{scriptsize}
\end{tikzpicture}
\caption{B\"or\"oczky configuration of $12$ lines}
\label{fig:B12}
\end{figure}
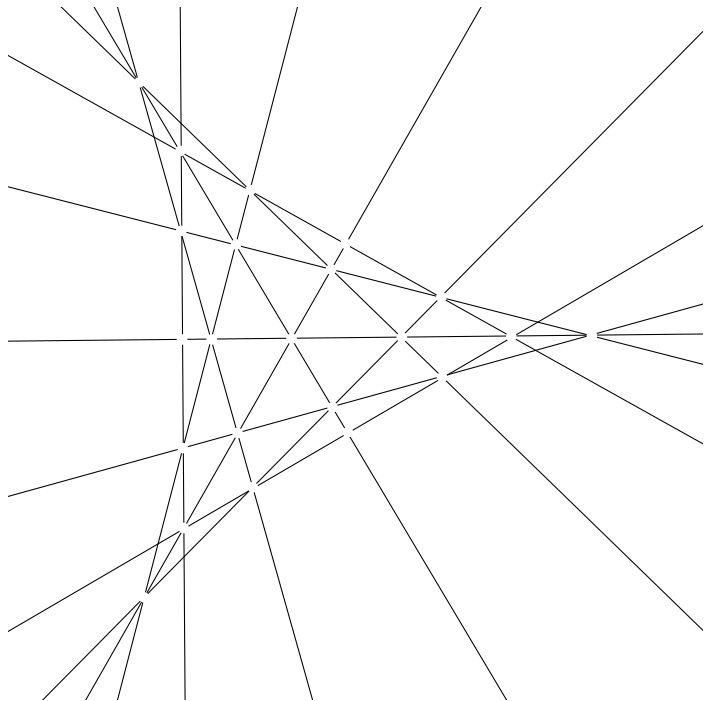
   Let $I$ be the saturated radical ideal defined by $Z_{12}$.
   Then
   $$I^{(3)}\nsubseteq I^2.$$
\end{cexample}
\proof
   The main idea is again to show that the product
   $$L_1\cdot L_2\cdot\ldots\cdot L_{9}$$
   of linear forms defining the lines $L_i$ is an element in $I^{(3)}$
   which is not contained in $I^2$. See \cite[Theorem 3]{Real}.
\endproof
   Contrary to the dual Hesse configuration, the configuration of $12$ lines outlined above
   allows much more freedom. It has been noticed already in \cite{DHNSST15}
   that it can be perturbed into a configuration of lines in the \emph{rational projective plane},
   see the discussion on page 389 of \cite{DHNSST15}. The moduli space of all such
   configurations has been studied in more detail by Lampa-Baczy\'nska and the second author with the result
   that they are parameterized by a three-dimensional rational variety,
   see \cite[Theorem A]{LBS15}. It came as a surprise that the analogous configuration
   of $15$ lines allows only $1$ parameter moving along an explicitly determined elliptic curve.
   In that case there are no rational configurations possible, see \cite[Theorem B]{LBS15}.

   Turning back to a more general setting, we observe first that
   Counterexample \ref{cex:dual hesse} is a special case of a whole series
   of configurations of lines leading to counterexamples to Problem \ref{pro:Huneke}.
   This series has been studied (with a different motivation) by Urzua in \cite[Example II.6]{Urz08}.
\begin{cexample}[Fermat configurations]\label{cex:Fermat}
   Let $n\geq 3$ be an integer and let
   $\K$ be a field of characteristic different from $2$ containing $n$ distinct roots
   of $1$. The polynomial
   $$F_n(x,y,z)=(x^n-y^n)(y^n-z^n)(z^n-x^n)$$
   splits completely into $3n$ linear forms. They define lines which
   intersect in $3$ points of multiplicity $n$ and $n^2$ points of multiplicity $3$.
   If $\eta$ is a primitive root of $1$ of order $n$,
   then the intersection points of these lines are
   $Q_{a,b}=(1:\eta^a:\eta^b)$ for $a,b=1\dots,n$
   and the three coordinate points $P_1=(1:0:0)$, $P_2=(0:1:0)$
   and $P_3=(0:0:1)$. There are exactly $n$ lines meeting
   in a coordinate point and exactly $3$ lines passing through
   every point $Q_{a,b}$. Taking $I$ as the radical ideal
   of the union of all points $Q_{a,b}$ and $P_i$, it is easy to see
   that $F_n\in I^{(3)}$. It is also easy to identify $I$ explicitly as
   $$I=\langle x(y^n-z^n), y(z^n-x^n), z(x^n-y^n)\rangle.$$
   It requires however some effort to show
   that $F_n$ is not contained in $I^2$ so that we have
   $$I^{(3)}\nsubseteq I^2.$$
\end{cexample}
\proof
   This is Proposition 2.1 in \cite{HarSec15}. An alternative proof is presented in
   \cite[Proposition 4.1]{Sec15}.
\endproof
   Additional counterexamples are in a sense isolated and they come from reflection groups
   acting on $\P^2$.
\begin{cexample}[Klein configuration]\label{cex:Klein}
   The configuration of bitangents of a smooth quartic curve in $\P^2$ has been a classical
   object of study. In particular for the Klein curve given by the equation
   $$xy^3+yz^3+zx^3=0$$
   it is highly symmetric. This has been first discovered and studied by Klein in \cite{Kle1879}.
   The Klein curve is the unique plane quartic curve with the
   group of automorphisms of the maximal size of $168$ elements, see \cite{JOS94}
   for a modern treatment.
   This group is the $\PSL_3(2)$ and it contains $21$ involutions whose fixed lines
   form a configuration of $21$ lines whose all incidences are $21$ points where 
   the lines meet by $4$ and $28$ points incident each to $3$ lines. This is the
   Klein configuration, one of dew known configurations with no simple
   intersection points, i.e. points incident to two lines only. 
   Taking $I$ to be the ideal of all $49$ intersection points mentioned above,
   we have
   $$I^{(3)}\nsubseteq I^2.$$
\end{cexample}
\proof
   See \cite[Theorem 4.4]{Sec15} for an even stronger statement.
\endproof
   The next configuration arises also from involutions in a reflection group of order $360$,
   see \cite{Wim1896}.
\begin{cexample}[Wiman configuration]
   The Wiman configuration consists of $45$ lines which intersect $120$ triple,
   $45$ quadruple and $36$ quintuple points. Taking as $I$ the ideal of all these
   intersection points we have
   $$I^{(3)}\nsubseteq I^2.$$
\end{cexample}
\proof
   This claim has been so far verified only by computer. We used Singular \cite{Singular}.
\endproof
\begin{remark}[A link to the Bounded Negativity Conjecture]
   Interestingly the last two counterexamples have appeared recently also
   in a seemingly unrelated topic of bounding Harbourne constants on
   birational models of $\P^2$, see \cite{BNAL}. It would be desirable
   to understand if there are some closer links between the Bounded Negativity Conjecture 
   and the Containment Problem.
\end{remark}
   Before we pass to counterexamples in finite characteristic we point out the following.
\begin{remark}\label{rmk:only r=2 m=3 N=2}
   In characteristic $0$ so far there are no counterexamples
   to Problem \ref{pro:BHH}, apart of those for $r=2$, $m=3$ and $N=2$.
\end{remark}
   The situation is quite different for points in projective planes in finite characteristic.
   Shortly after Counterexample \ref{cex:dual hesse} has been found, Bocci,
   Cooper and Harbourne pointed out in \cite[Remark 3.11]{BCH14} that there is 
   a similar example in characteristic $3$. Indeed the combinatorics of the two configurations
   is exactly the same. Nevertheless there are discrepancies between other
   algebraic invariant, e.g. the resurgence, see \cite[Theorem 2.1 and Theorem 3.2]{DHNSST15}.
\begin{cexample}[Lines in $\P^2(\F_3)$]
   Let $Z$ consist of all but one points in $\P^2(\F_3)$. Then there are
   $12$ points in $Z$, in each of which exactly $3$ out of $9$ lines missing
   the removed point meet. Then the containment \ref{eq:containment 3 in 2}
   fails for the ideal of $Z$.
\end{cexample}
\proof
   This is of course a special case of Counterexample \ref{cex:finite char in I2}.
\endproof
   The next counterexample is a finite characteristic version of Counterexample \ref{cex:Klein}.
\begin{cexample}[Conic in $\P^2(\F_7)$]
   It is well known that the Klein quartic reduces $\mod 7$ to a smooth (double)
   conic $C$ in $\P^2(F_7)$. The bitangents of the quartic correspond then to
   the secant lines of the conic. There are $28$ of them. There are $8$ more
   tangent lines to the conic and $21$ lines which do not meet $C$. These $21$ lines
   form again a configuration with $21$ quadruple intersection points and
   $28$ triple intersection points. Taking as $Z$ the set of all these
   intersection points the containment
   $$I^{(3)}\subset I^2$$
   fails for the ideal $I$ of $Z$.
\end{cexample}   
   
\section{Alternative approach}
   Problem \ref{pro:BHH} focuses on exponents on both sides of the containment
   $$I^{(m)}\subset I^r.$$
   It is however possible to sharpen the assertions of Theorem \ref{thm:ELS}
   by making the ideal on the right smaller in another way. This line of investigation
   was suggested in \cite{HaHu13}. This approach was motivated by the following
   conjecture due to Eisenbud and Mazur.
\begin{conjecture}[Eisenbud-Mazur]\label{conj:EM}
   Let $P\subset\C[[x_1,\ldots,x_N]]$ be a prime ideal in the ring
   of formal power series. Then
   $$P^{(2)}\subset M\cdot P,$$
   where $M=\langle x_1,\ldots,x_N\rangle$.
\end{conjecture}
   Thus the key idea in this path of investigation, paralleling
   Problem \ref{pro:containment} is to study the following question,
   see \cite[Question 1.3]{HaHu13}.
\begin{problem}[Harbourne, Huneke]\label{pro:containment HaHu}
   Decide for which $m$, $r$ and $j$ there is the containment
   $$I^{(m)}\subset M^jI^r.$$
\end{problem}
   For ideals defining fat points, Harbourne and Huneke formulated
   the following bold statement, see \cite[Conjecture 2.1]{HaHu13}.
\begin{conjecture}[Harbourne-Huneke]\label{conj:HaHu}
   Let $I\subset\K[x_0,\ldots,x_N]$ be a fat points ideal and let $r$ be a positive integer. Then the containment
   $$I^{(m)}\subset M^{r(N-1)}I^r$$
   holds for all $m\geq rN$.
\end{conjecture}
   Positive results towards this Conjecture have been obtained in
   \cite[Proposition 3.3]{HaHu13} for fat points ideals in $\P^2$
   arising as
   symbolic powers of radical ideals generated in a single degree.
   Particular examples of such ideals are star configurations of points
   and general points whose number is a binomial coefficient.

   Another classes of ideals supporting Conjecture \ref{conj:HaHu}
   have been studied in \cite[Theorem B]{DST13c}.

   To the best of our knowledge there is no counterexample
   to Conjecture \ref{conj:HaHu} known at present.
\section{A list of potential problems to think about}
   The article \cite{BCH14} contains a long list of conjectures, including
   Problem \ref{pro:BHH} and Conjecture \ref{conj:HaHu} discussed in the present note.
   Some of these conjectures require adjustments since the counterexamples
   presented here show their failure as well. These conjectures are also
   redundant in the sense that under some additional conditions some of them
   are stronger than the others. The picture presented in \cite{BCH14}
   is a little bit messy. For this reason, rather than repeating the conjectures
   and attempting to account on their current state, we have decided to present here
   a short list of open problems which are probably less challenging than
   Conjecture \ref{conj:HaHu} on the one hand but also more accessible
   on the other hand.

   It would be desirable to improve the Postulation Containment Criterion.
   The left hand side in the condition in Proposition \ref{prop:postulation containment criterion}
   is often used in order to bound the regularity of the power $I^r$.
   It is natural to wonder if one can use this regularity directly.
\begin{question}
   Can one replace $r\reg(I)$ by $\reg(I^r)$ in Proposition \ref{prop:postulation containment criterion}?
\end{question}
   The next question presents a similar approach based on another interesting invariant studied recently.
   Cutkosky, Herzog and Trung showed that
   the regularity of ordinary powers of a homogeneous ideal is eventually linear.
   An analogous question for asymptotic regularity of symbolic powers of ideals is
   far more involved. In the case relevant here, i.e. in the case of ideals with
   zero-dimensional support it has been studied recently by Cutkosky and Kurano in \cite{CutKur11}
   in a slightly more general setting of weighted projective spaces.
\begin{definition}[Symbolic asymptotic regularity]\label{def:symassreg}
   Let $I$ be a homogeneous ideal.
   The real number
   $$\symassreg(I)=\lim\limits_{m\to\infty}\frac{\reg(I^{(m)})}{m}$$
   is the \emph{symbolic asymptotic regularity} of $I$.
\end{definition}
   Cutkosky and Kurano showed that the limes in Definition \ref{def:symassreg}
   actually exists and it is equal to the so called $s$-invariant of $I$,
   which in turn is the reciprocal of the multi-point Seshadri constant of
   the reduced subscheme defined by $I$.
\begin{question}
   Can one replace $\reg(I)$ by $\symassreg(I)$ in Proposition \ref{prop:postulation containment criterion}?
\end{question}
   In the view of counterexamples presented in this note, an obvious attempt
   in order to save as much as possible of Problem \ref{pro:BHH} is the following statement.
\begin{question}
   Let $I$ be a homogeneous radical ideal of a finite set of points in $\P^N$ with $N\geq 3$.
   Let $r$ be a positive integer. Does then the containment
   $$I^{(m)}\subset I^r$$
   hold for all $m\geq Nr-1$?
\end{question}
   Thus a more specific question replacing in a sense Problem \ref{pro:Huneke}
   is the following.
\begin{question}
   Let $I$ be the radical ideal of a finite set of points in $\P^3$.
   Is then always
   $$I^{(5)}\subset I^2?$$
\end{question}
   The next problem is motivated by Conjecture \ref{conj:EM} even though it would
   not be a simple consequence of the Conjecture, see \cite[Question 2.5]{BCH14}.
\begin{question}
   Let $I$ be an arbitrary homogeneous ideal in $R$. Is then
   $$I^{(j+1)}\subset MI^{(j)}$$
   for all $j\geq 1$?
\end{question}
\paragraph*{Acknowledgement.}
   Research presented in this note has been started during the workshop
   ''Recent advances in Linear series and Newton-Okounkov bodies''
   held in Universit\'a degli Studi di Padova in February 2015.
   We would like to thank the organizers for invitation to participate
   in the workshop and for excellent working conditions and stimulating
   atmosphere they have created.
   We were partially supported by National Science Centre, Poland, grant
   2014/15/B/ST1/02197.


\bigskip \small

\bigskip
   Tomasz Szemberg,
   Polish Academy of Sciences, Institute of Mathematics, ul. \'Sniadeckich 8, PL-00-656 Warszawa, Poland

\nopagebreak
   \textit{E-mail address:} \texttt{tomasz.szemberg@gmail.com}

\bigskip
   Justyna Szpond,
   Department of Mathematics, Pedagogical University of Cracow,
   Podchor\c a\.zych 2,
   PL-30-084 Krak\'ow, Poland

\nopagebreak
   \textit{E-mail address:} \texttt{szpond@up.krakow.pl}


\end{document}